\documentclass{amsart}
\usepackage{epsfig}
\usepackage{amsmath}
\newtheorem{theorem}{Theorem}
\newtheorem{prop}[theorem]{Proposition}

\newtheorem{corollary}[theorem]{Corollary}
\newcommand\beq{\begin{equation}}
\newcommand\eeq{\end{equation}}
\newcommand\bce{\begin{center}}
\newcommand\ece{\end{center}}
\newcommand\bea{\begin{eqnarray}}
\newcommand\eea{\end{eqnarray}}
\newcommand\ben{\begin{enumerate}}
\newcommand\een{\end{enumerate}}

\newcommand\bs{\bigskip}
\newcommand\ms{\medskip}

\newcommand\brr{\begin{array}}
\newcommand\err{\end{array}}
\newcommand\bt{\begin{tabular}}
\newcommand\et{\end{tabular}}
\renewcommand\S{{\mathcal S}}
\newcommand\D{{\mathcal D}}
\newcommand\M{{\mathcal M}}
\newcommand\MM{{\mathcal N}}
\newcommand\T{{\mathcal T}}
\newcommand\bj{\Psi}
\newcommand\bjr{\Psi'}
\newcommand\bija{\varphi_1}
\newcommand\bijb{\varphi_2}
\newcommand\bijc{\varphi_3}
\newcommand\bijd{\varphi_4}
\newcommand\psiuc{\psi}
\newcommand\pre{\theta}
\newcommand\dgr{\delta}
\newcommand\krat{\Phi}

\title{Old and young leaves on plane trees}

\author{William Y. C. Chen}
\address{Center for Combinatorics, LPMC,
Nankai University, Tianjin 300071, P. R. China.}
\email{chen@nankai.edu.cn}
\author{Emeric Deutsch}
\address{Polytechnic University, Brooklyn, NY 11201.}
\email{deutsch@duke.poly.edu}
\author{Sergi Elizalde}
\address{Massachusetts Institute of Technology, Cambridge, MA 02139; Mathematical Sciences Research Institute, Berkeley, CA 94720.}
\email{elizalde@msri.org}

\begin{document}

\begin{abstract}
A leaf of a plane tree is called an old leaf if it is the leftmost child of its
parent, and it is called a young leaf otherwise. In this paper we enumerate
plane trees with a given number of old leaves and young leaves. The formula
is obtained combinatorially by presenting two bijections
between plane trees and 2-Motzkin paths which map young leaves to red
horizontal steps, and old leaves to up steps plus one. We derive
some implications to the enumeration of restricted permutations
with respect to certain statistics such as pairs of consecutive
deficiencies, double descents, and ascending runs. Finally, our main bijection
is applied to obtain refinements of two identities of Coker, involving refined
Narayana numbers and the Catalan numbers.
\end{abstract}

\maketitle

\section{Introduction}

Plane trees, also referred to as ordered trees, are a basic
object frequently used in combinatorics. Many enumerative results about them appear
throughout the literature. 
For example, a well-known interpretation of the Narayana numbers
is that they count the number of plane trees with a fixed number of
leaves. In this paper we classify the leaves of a plane tree into two different kinds,
distinguishing between old leaves and young leaves. This definition, which is
introduced in Section~\ref{sec:prel}, naturally gives rise to a refinement of the Narayana numbers.

These refined Narayana numbers also appear in the enumeration of 2-Motkin paths 
with respect to the number of up steps and red horizontal steps.
Such paths were introduced in \cite{BLPP},
and its structure has proved to be useful in the study of lattice paths,
noncrossing partitions, plane trees \cite{DeSh}, and other combinatorial 
objects and identities.
Our paper gives yet another example of the applicability of 2-Motzkin paths. The key to several of our results 
is a new bijection between plane trees and 2-Motzkin paths, with very convenient
properties. It provides a combinatorial derivation of the expression for the number of plane trees
with a given number of old and young leaves.

Another application of our bijection appears in \cite{CYY}, where Chen,
Yan and Yang use it to give combinatorial interpretations of two identities
involving the Narayana numbers and Catalan numbers, due to Coker \cite{Cok}.
This way they solve the two open problems left in \cite{Cok}. Here we will show that 
a more detailed analysis of the bijection and its properties 
gives refinements of the two identities of Coker, as well as bijective
proofs of these refinements.

The paper is structured as follows. In Section~\ref{sec:prel} we review some definitions and notation
about plane trees, Dyck paths, Motzkin paths, and 2-Motzkin paths.
We also introduce the concepts of old leaves and young leaves of a
plane tree. In Section~\ref{sec:enum} we give the generating
function for plane trees with variables marking the number of old
leaves and the number of young leaves, as well as exact formulas for
the number of plane trees of a given size when the number of old and
young leaves is fixed. In Section~\ref{sec:bj} we present two
bijections from the set of plane trees with $n$ edges to the set of
2-Motzkin paths of length $n-1$. Some interesting properties  of
these bijections are studied in Section~\ref{sec:conseq}. We show
that they map old and young leaves of trees into statistics on
2-Motzkin paths that are easier to deal with. In
Section~\ref{sec:perms} we describe some bijections between plane
trees and permutations avoiding patterns of length 3, and
investigate what old and young leaves are mapped to by these
bijections. This implies that the distribution of certain parameters
on restricted permutations is given by the same formulas
enumerating plane trees with respect to old and young leaves.
Finally, in Section~\ref{sec:coker} we apply our bijection to obtain
refinements of two combinatorial identities due to Coker \cite{Cok} and
proven combinatorially by Chen, Yan and Yang~\cite{CYY}.

\section{Preliminaries}\label{sec:prel}

\subsection{Plane trees}

A \emph{plane tree} $T$ can be defined recursively as a finite set
of vertices such that one distinguished vertex $r$ is called the
\emph{root} of $T$, and the remaining vertices are put into an
ordered partition $(T_1,T_2,\ldots,T_m)$ of $m$ disjoint non-empty
sets, each of which is a plane tree. We will draw plane trees with
the root on the top level, with edges connecting it to the roots of
$T_1,T_2,\ldots,T_m$, which will be drawn from left to right on the
second level. For each vertex $v$, the nodes in the next lower level
connected to $v$ by an edge are called the \emph{children} or
\emph{successors} of $v$, and $v$ is called the \emph{parent} of its
children. Clearly each vertex other than $r$ has exactly one parent.
A vertex of $T$ is called a \emph{leaf} if it has no children (by
convention, we assume that the \emph{empty} tree, formed by a single
node, has no leaves).

We denote by $\T_n$ the set of plane trees with $n$ edges.
It is well-known that $|\T_n|=\frac{1}{n+1}\binom{2n}{n}$, the $n$-th Catalan
number, and that the number of trees with $n$ edges and $k$ leaves
is the Narayana number $\frac{1}{n}\binom{n}{k}\binom{n}{k-1}$.

\begin{figure}[hbt]
\epsfig{file=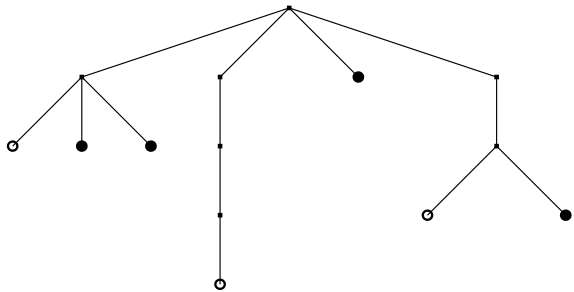,width=3in} \caption{\label{fig:tree}
A tree with 3 old leaves and 4 young leaves.}
\end{figure}

We classify the leaves of a plane tree into old and young leaves.
We say that a leaf is an \emph{old leaf} if it is the leftmost
child of its parent, and that it is a \emph{young leaf} otherwise.
For example, the tree in Figure~\ref{fig:tree} has four young leaves,
drawn with black filled circles, and three old leaves, drawn with empty circles.
The enumeration of plane trees with respect to the number of old
and young leaves is done in Section~\ref{sec:enum}.

\subsection{Lattice paths}

We review the definitions of Dyck, Motzkin, and 2-Motzkin paths.
They are all lattice paths in $\mathbb{Z}^2$ starting at $(0,0)$,
ending on the $x$-axis, and never going below this axis. A
\emph{Dyck path} consists of steps $U=(1,1)$ 
and $D=(1,-1)$. 
In a \emph{Motzkin path} we allow also horizontal steps $H=(1,0)$,
so that the path is a sequence of steps $U$, $D$ and $H$.
A \emph{2-Motzkin path} consists of up
and down steps, and horizontal steps that can be colored either red
or blue. We use $R$ to denote a red step, and $B$ for a blue step.
In the pictures in this paper, red steps will be drawn with a dashed
line to make them clearly distinguishable from blue steps, which will be
drawn with a solid line.
The length of any of these paths is the total number of steps.

We shall denote by $\D_n$ the set of Dyck paths of length $2n$, by
$\M_n$ the set of Motzkin paths of length $n$, and by $\MM_n$ the
set of 2-Motzkin paths of length $n$. The number of paths of each
kind is given by $|\D_n|=C_n$, $|\M_n|=M_n$, and $|\MM_n|=C_{n+1}$,
where $M_n=\sum_{k=0}^n \binom{n}{2k} C_k$ is the $n$-th Motzkin
number.

The generating function for Catalan numbers is
$C(z)=\sum_{n\ge0}C_nz^n=\frac{1-\sqrt{1-4z}}{2z}$, and the one
for Motzkin numbers is $M(z)=\sum_{n\ge0}M_nz^n=\frac{1 - z -
\sqrt{1-2z-3z^2}}{2z^2}$.

\section{Enumeration of trees with respect to old and young leaves} \label{sec:enum}

Here we give an expression for the generating function
$$G(t,s,z)=\sum_T t^{\mathrm{\# old\,leaves\,of\,T}}
s^{\mathrm{\# young\,leaves\,of\,T}} z^{\mathrm{\#
edges\,of\,T}},$$ where the sum is over all plane trees $T$, and
$t$ and $s$ mark the number of old and young leaves respectively.

\begin{theorem} Let $G(t,s,z)$ be defined as above. We have
$$G(t,s,z)=\frac{1+z-sz-\sqrt{1-2(1+s)z+(1-4t+2s+s^2)z^2}}{2z}.$$
\end{theorem}

\begin{proof}
We will find an equation for $G(t,s,z)$ using a decomposition of
plane trees. Let $T$ be any plane tree, and let $m$ be the number of
children of the root. If $m=0$, then the tree has no edges, and its
contribution to the generating function $G$ is $1$. If $m\ge1$, let
$T_1,T_2,\ldots,T_m$ be the sequence of plane trees hanging from
left to right from the children of the root. If $T_1$ has no edges,
then it creates an old leaf of $T$, otherwise all the old (resp.
young) leaves of $T_1$ become old (resp. young) leaves of $T$.
Therefore, the contribution to the generating function of $T_1$ and
the edge connecting it to the root is $z(G(t,s,z)-1+t)$. For
$i\ge2$, old and young leaves of $T_i$ become leaves of $T$ of the
same kind as well. However, if $T_i$ is has no edges, then an
additional young leaf of $T$ is created. Thus, the contribution to
the generating function of each $T_i$ with $i\ge2$ and the edge
connecting it to the root is $z(G(t,s,z)-1+s)$. It follows that for
$m\ge1$, the contribution of the plane trees whose root has degree
$m$ is $z^m(G-1+t)(G-1+s)^{m-1}$. Summing over all $m\ge0$ we obtain
\beq\label{eq:G}
G(t,s,z)=1+\frac{z(G(t,s,z)-1+t)}{1-z(G(t,s,z)-1+s)}.\eeq Isolating
$G$ the formula follows.
\end{proof}

\begin{prop}\label{prop:numbers} \ben
\item The number of plane trees with $n$ edges, $i$ old leaves, and
$j$ young leaves is
$$\frac{1}{n}\binom{n}{i}\binom{n-i}{j}\binom{n-i-j}{i-1}.$$
\item The number of plane trees with $n$ edges and $k$ old leaves
is $$\frac{2^{n-2k+1}}{k}\binom{n-1}{2k-2}\binom{2k-2}{k-1}.$$
\item The number of plane trees with $n$ edges and $k$ young leaves
is $$\binom{n-1}{k}M_{n-k-1}.$$ \een
\end{prop}

\begin{proof}[First proof]
Applying Lagrange inversion formula to equation (\ref{eq:G}), we
obtain that the coefficient of $t^is^jz^n$ in $G(t,s,z)$ is
$\frac{1}{n}\binom{n}{i}\binom{n-i}{j}\binom{n-i-j}{i-1}$, which
is the first expression. For the other two expressions, apply
Lagrange inversion to the same equation, after the substitutions
$t=1$ and $s=1$ respectively.
\end{proof}

\begin{proof}[Second proof]
We can give a bijective proof of the first part of Proposition~\ref{prop:numbers} as follows.
In \cite{Ch90}, the author gives a bijective algorithm to decompose
any labeled plane tree with $n$ edges into a set $F$ of $n$ matches
with labels $\{1, \ldots, n, n+1, (n+2)^*, \ldots, (2n)^*\}$, where
a match is a rooted tree with two vertices. The reverse procedure of
the decomposition algorithm is the following merging algorithm. We
start with a set $F$ of matches on $\{1, \ldots, n+1, (n+2)^*,
\ldots, (2n)^*\}$. A vertex labeled by a mark $*$ is called a marked
vertex.

\noindent (1) Find the tree $T$ with the smallest root in which no
vertex is marked. Let $i$ be the root of
$T$.\\
(2) Find the tree $T^*$ in $F$ that contains the smallest marked
vertex. Let $j^*$ be this marked vertex.\\
(3) If $j^*$ is the root of $T^*$, then merge $T$ and $T^*$ by
identifying $i$ and $j^*$, keep $i$ as the new vertex, and place the
subtrees of $T^*$ to the right of $T$. The operation is called a
\emph{horizontal merge}. If $j^*$ is a leaf of $T^*$, then replace
$j^*$ with $T$ in $T^*$. This operation is called a \emph{vertical
merge}. See Figure \ref{merge}.
\begin{figure}[h,t]
\begin{center}
\begin{picture}(420,50)
\setlength{\unitlength}{1.1mm} \linethickness{0.4pt}

\put(5,10){\line(-1,-1){5}} \put(0,5){\circle*{1}}
\put(0,5){\line(0,-1){5}} \put(0,0){\circle*{1}}
\put(5,10){\circle*{1}} \put(5,5){\circle*{1}}
\put(5,0){\line(0,1){10}} \put(5,0){\circle*{1}} \put(3,13){$i$}

\put(11,4){\large\bf$+$} \put(20,0){\circle*{1}}
\put(20,10){\circle*{1}} \put(20,0){\line(0,1){10}}
\put(18,13){$j^*$} \put(25,4){\large\bf$\Rightarrow$}

\put(35,0){\circle*{1}} \put(35,5){\circle*{1}}
\put(35,0){\line(0,1){5}} \put(40,10){\circle*{1}}
\put(35,5){\line(1,1){5}} \put(40,5){\circle*{1}}
\put(40,10){\line(0,-1){10}} \put(40,0){\circle*{1}}
\put(40,10){\line(1,-1){5}} \put(45,5){\circle*{1}}

\put(39,13){$i$} \put(2,-7){$T$} \put(18,-7){$T^*$}

\put(80,0){\circle*{1}} \put(80,0){\line(1,2){5}}
\put(85,10){\circle*{1}} \put(90,0){\line(-1,2){5}}
\put(90,0){\circle*{1}} \put(83,13){$i$}

\put(95,4){\large\bf$+$}

\put(105,0){\circle*{1}} \put(105,10){\circle*{1}}
\put(105,0){\line(0,1){10}} \put(106,-1){$j^*$}
\put(110,4){\large\bf$\Rightarrow$} \put(120,0){\circle*{1}}
\put(125,5){\circle*{1}} \put(120,0){\line(1,1){5}}
\put(130,0){\circle*{1}} \put(130,0){\line(-1,1){5}}
\put(125,10){\circle*{1}} \put(125,5){\line(0,1){5}}
\put(127,4){$i$} \put(85,-7){$T$} \put(103,-7){$T^*$}
\end{picture}
\end{center}
\ms
\caption{Horizontal merge and vertical merge}\label{merge}
\end{figure}
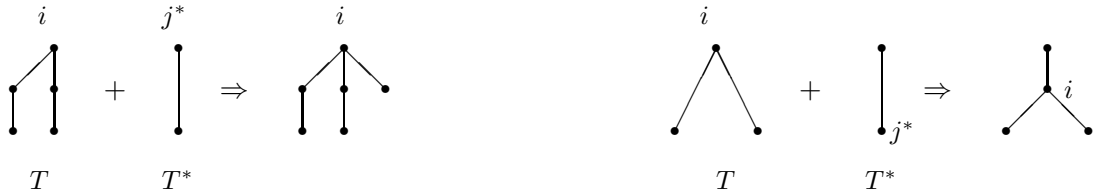

\noindent (4) Repeat the above procedure until $F$ becomes a
labeled tree.

For any labeled plane tree with $n$ edges, $i$ old leaves, and
$j$ young leaves, the corresponding set $F$ of $n$ matches
consists of $i$ matches labeled by labels without marked
vertices,
 $j$ matches with marked roots and unmarked leaves, and all leaves in
the remaining matches are  marked vertices. Thus, we can count the
number of labeled plane trees with $n$ edges, $i$ old leaves, and
$j$ young leaves as follows:
\begin{eqnarray*}
{n+1 \choose 2i} \frac{(2i)!}{i!} {n+1-2i \choose j}{n-1 \choose
j}j! {n-1-j \choose n-i-j}(n-i-j)! = \frac{(n+1)!}{n}{n \choose
i}{n-i \choose j}{n-i-j \choose i-1}.
\end{eqnarray*}
Now, to count unlabeled plane trees we just divide by $(n+1)!$, obtaining the
desired formula.

Summing for all $j$ we obtain the formula in part (2) of the proposition, and
summing for all $i$ we derive the third formula.
\end{proof}

Particular cases of this proposition give rise to the following
two statements. The second one appeared already in~\cite{DS77} as
a manifestation of the Motzkin numbers.

\begin{corollary}\label{cor:noyoung}\ben
\item The number of plane trees in $\T_n$ with exactly one old
leaf is $2^{n-1}$.
\item The number of plane trees in $\T_n$ with no young leaves
is $M_{n-1}$. \een
\end{corollary}

\section{Two bijections between plane trees and 2-Motzkin paths} \label{sec:bj}

In this section we present two bijections $\bj$ and $\bjr$ from the set of plane
trees with $n$ edges and the set of 2-Motzkin paths of length
$n-1$.
These bijections have the convenient property that they map old
and young leaves of the tree to certain statistics of the
2-Motzkin path that are very easy to deal with, as shown in the
next section. This will allow us to give bijective proofs of
Corollary~\ref{cor:noyoung} and some parts of
Proposition~\ref{prop:numbers}.
Both bijections have very similar properties, and in fact one of them would be
enough to prove the results in the next section. However, they are
defined in quite different ways, and we feel that presenting both
bijections gives a better insight on how old and young leaves correspond
to statistics on paths.

\ms

Let us start describing the bijection $\bj$. It consists of three
steps. Given a plane tree $T\in\T_n$ (assume $n\ge1$), we first
transform it into a Dyck path using the following well-known
bijection, which we denote $\pre$. Starting from the root, traverse
the edges of the tree in preorder from right to left. To each edge
passed on the way down there corresponds a step $U$, and to each
edge passed on the way up there corresponds a step $D$. This gives
us a Dyck path $\pre(T)$ of length $2n$.

The next step is to replace each peak $UD$ of the path followed by
a $U$ step with a red horizontal step $R$. That is, we traverse
the path $\pre(T)$ from left to right replacing each $UDU$ with $RU$. This
gives us a Motzkin path with steps $U$, $D$ and $R$, whose length
is variable.

Finally, we need to transform this Motzkin path into a 2-Motzkin
path $\bj(T)$ of length $n-1$. The bijection that we will use for
this purpose is essentially the same one described by
Callan~\cite{Ca04} between $UDU$-free Dyck paths and Motzkin paths,
where we ``ignore" the steps $R$ of our path and let the new level
steps be all $B$ steps. Notice that after the transformation in the
previous paragraph, every peak $UD$ in our Motzkin path is followed
by a $D$ step, unless it is at the end of the path. This last
transformation is done as follows. Place a mark on each $U$ that is
followed by a $D$, on each $D$ that is followed by another $D$, and
on the $D$ at the end of the path. Next, change each unmarked $U$
whose matching $D$ is marked into an $B$. (The matching $D$ is the
step that is encountered directly east of $U$.) Lastly, delete all
the marked steps.

\begin{figure}[hbt]
\epsfig{file=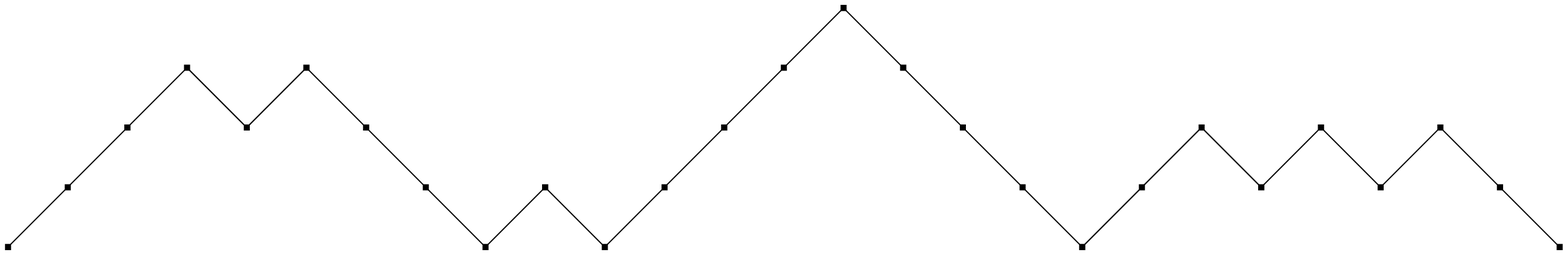,width=5.5in} \caption{\label{fig:dyck}
The Dyck path $\pre(T)$ for $T$ in Figure~\ref{fig:tree}.}
\end{figure}

\begin{figure}[hbt]
\epsfig{file=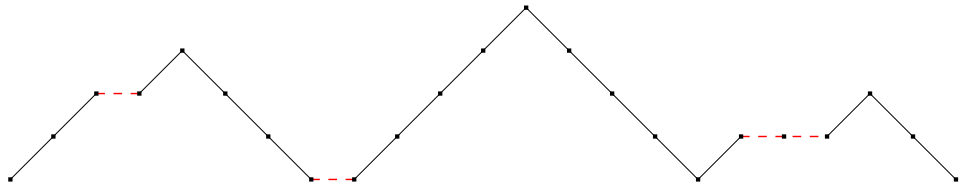,width=5in} \caption{\label{fig:motzkin}
The Motzkin path $UURUDDDRUUUUDDDDURRUDD$.}
\end{figure}

\begin{figure}[hbt]
\epsfig{file=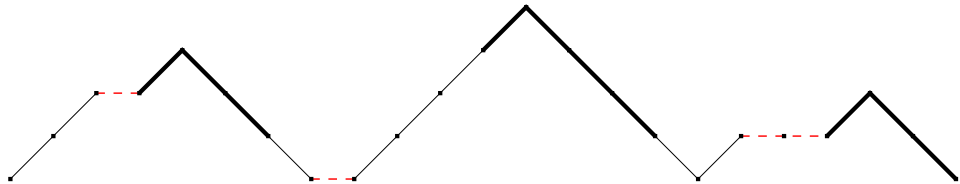,width=5in} \caption{\label{fig:marked}
The Motzkin path with some steps marked.}
\end{figure}

\begin{figure}[hbt]
\epsfig{file=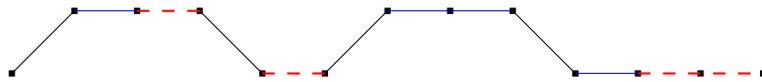,width=4in} \caption{\label{fig:2motzkin}
The 2-Motzkin path $\bj(T)=UBRDRUBBDBRR$.}
\end{figure}

After this procedure we obtain a 2-Motzkin path $\bj(T)$ with
$n-1$ steps. For example, for the tree $T$ in Figure~\ref{fig:tree},
applying the first part of the bijection we get the Dyck path
in Figure~\ref{fig:dyck}.
Replacing each $UDU$ with $RU$,
we get the Motzkin path in Figure~\ref{fig:motzkin}. In the third
part of the bijection, we mark the steps that in Figure~\ref{fig:marked} are thicker.
Changing each unmarked $U$ with a marked matching $D$ to a $B$, we
get
$UBR\dot{U}\dot{D}\dot{D}DRUBB\dot{U}\dot{D}\dot{D}\dot{D}DBRR\dot{U}\dot{D}\dot{D}$,
where the dots indicate the marked steps. Finally, deleting the
marked steps, we obtain the 2-Motzkin path in
Figure~\ref{fig:2motzkin}.

It is clear that the first two steps of this map are reversible,
that is, from the Motzkin path with steps $U$, $D$ and $R$ it is
easy to recover the tree. The fact that the last step is a
bijection as well follows from the description of the inverse
given in~\cite{Ca04}. The only difference here is that we need to
disregard the steps $R$ that we have now in the path, since they
are not affected by this part of the bijection.

\ms

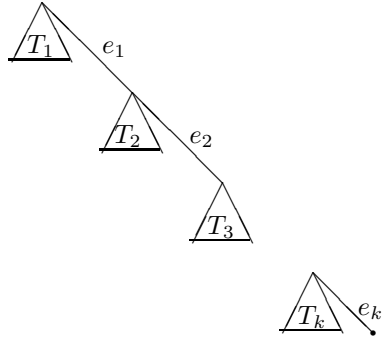
\begin{figure}[h,t]
\begin{center}
\begin{picture}(200,150)
\setlength{\unitlength}{0.4cm} \put(14,1){\circle*{0.2}}
 \put(12,3){\line(1,-1){2}}
\put(12,3){\line(-1,-2){1}} \put(12,3){\line(1,-2){1}}
\put(10.9,1.1){\line(1,0){2}} \put(11.5,1.3){$T_k$}
\put(13.5,1.6){$e_k$}
\multiput(12,3)(-0.4,0.4){7}{\line(-1,1){0.1}}
\put(9,6){\line(-1,-2){1}} \put(9,6){\line(1,-2){1}}
\put(7.9,4.1){\line(1,0){2}} \put(8.5,4.3){$T_{3}$}
\put(9,6){\line(-1,1){3}} \put(6,9){\line(-1,-2){1}}
\put(6,9){\line(1,-2){1}} \put(4.9,7.1){\line(1,0){2}}
\put(5.4,7.3){$T_2$} \put(7.9,7.3){$e_2$}
\put(6,9){\line(-1,1){3}} \put(3,12){\line(-1,-2){1}}
\put(3,12){\line(1,-2){1}} \put(1.9,10.1){\line(1,0){2}}
\put(2.5,10.3){$T_1$} \put(5,10.3){$e_1$}
\end{picture}
\end{center}
\caption{\label{fig:decomptree} Decomposition of a plane tree.}
\end{figure}

Now we describe another bijection $\bjr$ between $\T_n$ and the set of 2-Motzkin paths of length $n-1$.
We can construct $\bjr$ recursively. Given a plane tree $T$, consider the decomposition given in Figure~\ref{fig:decomptree},
where $e_1,e_2,\cdots,e_k$ is the path obtained starting at the root and successively descending to the rightmost child until we
reach a leaf. $T_1,T_2,\cdots,T_k$ are possibly empty subtrees hanging from the vertices of this path. When $i\neq k$,
if the subtree $T_i$ consists of a single vertex, then we encode
the edge $e_i$ by $B$; otherwise, $e_i$ is encoded by a $U$ and a $D$.
When $i=k$, if the subtree $T_i$ consists of a single vertex, then
the edge $e_k$ does not produce any step in the encoding; otherwise, $e_k$ is encoded
by $R$. We traverse the path from $e_1$ to $ e_k$ and construct $\bjr(T)$ as follows. If the encoding of
$e_i$ is $B$ or $R$, then we record $Q_i=B$ or
$Q_i=R\bjr(T_i)$ respectively. If the encoding of $e_i$ is a $U$ and a $D$, then we
denote $Q_i=U\bjr(T_i)D$. Joining these segments $Q_i$ from $1$ to
$k$, we obtain a $2$-Motzkin path $M=Q_1Q_2\cdots Q_k$.

Here is an alternative way to describe $\bjr$. Given a plane tree $T$ with $n$ edges, label its
vertices by $U, D, B$ or $R$ while traversing it in
preorder. For an internal vertex, if it is not the leftmost
child of its parent we label it by $U$, otherwise we label the
vertex by $B$. A young leaf is labeled by $R$ and an old leaf
by $D$, except the last old leaf that we encounter, which is left unlabeled.
This way all vertices get a label except the root and the last old leaf.

To construct now the $2$-Motzkin path we traverse the vertices of the tree in a different order and read the labels.
Suppose that the root of $T$ has $k$ children $v_1, v_2,
\ldots, v_k$ and that $T_i$ is the subtree with root $v_i$. Then we traverse first the vertices
$v_k, v_{k-1}, \ldots, v_1$ in this order, and then traverse $T_1,
T_2, \ldots, T_k$ recursively. It can be shown that the path obtained in this way is $\bjr(T)$.

\section{Consequences of the bijections} \label{sec:conseq}

The main properties of $\bj$ and $\bjr$ are given in the following
proposition. We state it only for $\bj$, but exactly the same result
holds if we replace $\bj$ with $\bjr$. The proof for $\bjr$ follows
easily from its recursive description.

\begin{prop}\label{prop:corresp} Let $T$ be a plane tree with $n\ge1$ edges, and let
$P=\bj(T)$ be the corresponding 2-Motzkin path. We have \ben
\item $\#$ of old leaves of $T = 1 + \#$ of $U$ steps of $P$,
\item $\#$ of young leaves of $T = \#$ of $R$ steps of $P$.
\een
\end{prop}

\begin{proof}
Let us first take a look at how old and young leaves are
transformed by the first part $\pre$ of the bijection, which consists in
reading $T$ in preorder from right to left and building a Dyck
path out of it. It is clear that each leaf of $T$ produces a peak
in $\pre(T)$. Now, a young leaf of $T$ corresponds to a peak $UD$
followed by a $U$ step, whereas an old leaf of $T$ gives rise to a
peak $UD$ not followed by a~$U$.

The second part of the bijection transforms each peak $UD$ followed
by a $U$ into a red step $R$, and these steps remain unchanged by
the third part of the bijection. This proves~(2). The remaining
peaks of the Dyck path are followed either by a $D$ or by nothing,
and they are not affected by the second part of the bijection, so
these are the only peaks in the Motzkin path. In the final part, we
place a mark on each $D$ that is followed by another $D$ or by
nothing, and the only $D$'s that are not erased are the unmarked
ones. Therefore, the number of $U$ steps (equivalently, the number
of $D$ steps) in $\bj(T)$ equals the number of $D$'s in the Motzkin
path that are left unmarked. The $D$ steps in the Motzkin path can
be grouped in sequences of consecutive $D$'s, each such sequence
immediately following a peak (note that the path has no occurrences
of $RD$, so each $D$ is in one of these sequences). In the sequence
of $D$'s following the rightmost peak all the steps are marked. For
each remaining peak, among the $D$ steps in the consecutive sequence
following it, all but the last one are marked. Thus, only one $D$
step survives for each peak other than the rightmost one. In other
words, the number of $D$ steps in $\bj(T)$ is the number of peaks of
the Motzkin path minus one. This implies~(1).
\end{proof}

By means of the bijection $\bj$ and the properties described above, we
can now give a combinatorial proof of Corollary~\ref{cor:noyoung}.
To prove the first part, observe that by property (1) of
Proposition~\ref{prop:corresp}, $\bj$ induces a bijection between
plane trees with exactly one old leaf and 2-Motzkin paths with no
$U$ steps. But these paths are just sequences of horizontal steps,
each of which can be colored red or blue. Thus, the number of
plane trees on $n$ edges with exactly one old leaf is $2^{n-1}$.

A direct proof of this nice fact, without using bijections to
lattice paths, can be given as follows. Let $T$ be a tree with $n$
edges and exactly one old leaf, call it $\ell$. We can find $\ell$
by following the path that starts at the root and always continues
to the leftmost child. Let $P$ be this path. Then $\ell$ must be
at the end of $P$. Now we claim that the remaining nodes of $T$
are leaves hanging from the nodes of $P$ other than $\ell$.
Indeed, if a node of $P$ had a child not in $P$ with successors,
then following the path that starts at this child and continues
always to the leftmost child, we would end at another old leaf,
which is a contradiction. Reciprocally, if only leaves are hanging
from $P$, then no more old leaves appear. Now, the number of trees
consisting of a path $P$ with leaves hanging from its nodes is
clearly $2^{n-1}$. Indeed, one can think of it as a composition of
$n$, say $n = a_1 + a_2 +\cdots$, where $a_i$ is the number of
children of the $i$-th node of $P$.

More generally, we can use our bijection to give a combinatorial
proof of the second part of Proposition~\ref{prop:numbers}, namely
that the number of plane trees with $n$ edges and $k$ old leaves
is $\frac{2^{n-2k+1}}{k}\binom{n-1}{2k-2}\binom{2k-2}{k-1}$. By
the first property of $\bj$ given above, we have to count the
number of 2-Motzkin paths of length $n-1$ with $k-1$ $U$ steps. To
produce such a path, we can choose in $\binom{n-1}{2k-2}$ ways the
positions of the $k-1$ $U$'s and $k-1$ $D$'s in the path. Deciding
which of these positions will be filled with a $U$ or with a $D$
is equivalent to choosing a Dyck path with $2k-2$ steps, and this
can be done in $\frac{1}{k}\binom{2k-2}{k-1}$ ways. The remaining
$n-2k+1$ positions are horizontal steps, which can be colored red
or blue in $2^{n-2k+1}$ ways.

To show the second part of Corollary~\ref{cor:noyoung}
combinatorially, notice that property (2) of
Proposition~\ref{prop:corresp} implies that $\bj$ maps plane trees
with no young leaves into 2-Motzkin paths with no $R$ steps. These
are just Motzkin paths with steps $U$, $D$ and $B$. Therefore, the
number of plane trees on $n$ edges with no young leaves equals the
number of Motzkin paths with $n-1$ steps, which is $M_{n-1}$.

More generally, the same property of $\bj$ can be used to prove
the last part of Proposition~\ref{prop:numbers}, namely that the
number of plane trees with $n$ edges and $k$ young leaves is
$\binom{n-1}{k}M_{n-k-1}$. Indeed, now the problem is equivalent
to counting 2-Motzkin paths of length $n-1$ with $k$ $R$ steps. We
can choose in $\binom{n-1}{k}$ ways where these $R$ steps go, and
then the remaining $n-k-1$ steps can be filled with a Motzkin path
with steps $U$, $D$ and $B$. 

{\bf Remark.} Another combinatorial proof of part (3) of 
Proposition~\ref{prop:numbers} can be obtained
using the result mentioned in \cite{DeSh} (and proved also in~\cite{Sun04})
that $\binom{n-1}{k}M_{n-k-1}$ counts the number of Dyck paths of
length $2n$ with $k$ $DUD$'s.

\ms

The description of $\bj$ implicitly contains a bijection between Dyck paths and 2-Motzkin paths.
There is a simpler bijection, perhaps the most standard one,
that transforms a 2-Motzkin path of length $n-1$ into a Dyck path of length $2n$, by first applying the following rules:
$$U\rightarrow UU,\quad D\rightarrow DD,\quad R\rightarrow UD,\quad B\rightarrow DU,$$
and then inserting a $U$ at the beginning and a $D$ at the end of the path.
Applying $\bj$ followed by this bijection, young leaves of the tree are mapped to peaks at even height in the Dyck path.
This shows that the statistic `number of young leaves' in $\T_n$ is equidistributed with the statistic `number of peaks at even height' in $\D_n$.

\section{Some statistics on restricted permutations} \label{sec:perms}

Using some known bijections between Dyck paths and permutations avoiding a pattern of length 3,
the parameters counting the number of old and young leaves in plane trees correspond to certain statistics on
restricted permutations. Given a pattern $\sigma$, we denote by $\S_n(\sigma)$ the set of
permutations in the symmetric group $\S_n$ avoiding $\sigma$. It is well-known \cite{Knu73} that
if $\sigma$ is any pattern of length 3, then $|\S_n(\sigma)|=C_n$, the $n$-th Catalan number.

We begin with a few definitions. Let $\pi$ be a permutation. We say that $\pi_i$ is an \emph{excedance}
if $\pi_i>i$, that it is a \emph{weak excedance} if $\pi_i\ge i$, and that it is a \emph{deficiency} if $\pi_i<i$..
A \emph{left-to-right minimum} of $\pi$ is an element $\pi_i$ such
that $\pi_i<\pi_j$ for all $j<i$. We call a \emph{double descent} of
$\pi$ a sequence of three consecutive decreasing elements
$\pi_i>\pi_{i+1}>\pi_{i+2}$ (equivalently, two consecutive
descents). A \emph{double ascent} is defined analogously. An
\emph{ascending run} is a maximal increasing sequence of (at least
two) consecutive elements of $\pi$, i.e.,
$\pi_i<\pi_{i+1}<\cdots<\pi_{i+k}$, with $k\ge1$.

\begin{prop}\label{prop:st321}
There is a bijection $\bija:\T_n\longrightarrow\S_n(321)$ such that, if $T\in\T_n$
and $\pi:=\bija(T)\in\S_n(321)$, then
\ben
\item $\#$ of young leaves of $T = \#$ of pairs of consecutive weak excedances of $\pi$,
\item $\#$ of old leaves of $T = \#$ of weak excedances of $\pi$ not followed by another weak excedance.
\een
\end{prop}

\begin{proof}
We use a bijection $\psiuc$ between $\S_n(321)$ and $\D_n$ which is very similar to the one
given by Krattenthaler \cite{Kra01} from $\S_n(123)$ to $\D_n$. Here is a way to describe it.
Let $\pi\in\S_n(321)$, and let $\pi_{i_1},\pi_{i_2},\ldots,\pi_{i_k}$ be its weak excedances, from left to right.
Define $\psiuc(\pi)$ to be the path that starts with $\pi_{i_1}$ up steps, then has, for each $j$
from 2 to $k$, $i_j-i_{j-1}$ down steps followed by
$\pi_{i_j}-\pi_{i_{j-1}}$ up steps, and finally ends with $n+1-i_k$
down steps. It can be checked that this is indeed a bijection between $321$-avoiding permutations and Dyck paths.

Our bijection $\bija$ is defined as  $\bija=\psiuc^{-1}\circ\pre$. Recall that $\pre$ reads a plane
tree in preorder from right to left and creates a Dyck path.

We saw that young leaves of $T$ correspond to occurrences of $UDU$ in the path $\pre(T)$, and that
old leaves of $T$ are mapped by $\pre$ to either a $UDD$ or a terminal (i.e., at the end of the path) $UD$.
Now, if $\pi\in\S_n(321)$, a $UDU$ is obtained in $\psiuc(\pi)$ precisely when we have a weak excedance
followed by another weak excedance, which causes one of the descending slopes to have length $i_j-i_{j-1}=1$.
Similarly, a $UDD$ corresponds to a weak excedance followed by a deficiency (i.e., an element that is not a weak
excedance), and a terminal $UD$ corresponds to the weak excedance $\pi_n=n$.
\end{proof}

For example, if $T$ is the tree in Figure~\ref{fig:tree}, with $\pre(T)$ given in Figure~\ref{fig:dyck}, then the corresponding
permutation is $\bija(T)=(3,4,1,2,5,9,6,7,8,11,12,13,10)\in\S_{12}(321)$. It has four pairs of consecutive weak excedances, namely $(3,4)$,
$(5,9)$, $(11,12)$ and $(12,13)$, and three weak excedances not followed by another weak excedance, namely $4$, $9$ and $13$.

A similar result for $132$-avoiding permutations is given next. For $\pi\in\S_n$, let $(n+1)\pi$ (resp. $\pi(n+1)$) be the permutation in $\S_{n+1}$
obtained by inserting $n+1$ at the beginning (resp. at the end) of $\pi$.

\begin{prop}
There is a bijection $\bijb:\T_n\longrightarrow\S_n(132)$ such that, if $T\in\T_n$
and $\pi:=\bijb(T)\in\S_n(132)$, then
\ben
\item $\#$ of young leaves of $T = \#$ of double descents of $(n+1)\pi$,
\item $\#$ of old leaves of $T = \#$ of ascending runs of $\pi(n+1)$.
\een
\end{prop}

\begin{proof}\label{prop:st132}
We use the bijection from $\S_n(132)$ to $\D_n$ denoted by $\krat$ that appears in Krattenthaler \cite{Kra01}.
Given $\pi\in\S_n(132)$, let $\pi_{i_1},\pi_{i_2},\ldots,\pi_{i_k}$ be its left-to-right minima, from left to right.
Then $\krat(\pi)$ is the Dyck path that starts with $n+1-\pi_{i_1}$ up steps, then has, for each $j$
from 2 to $k$, $i_j-i_{j-1}$ down steps followed by
$\pi_{i_{j-1}}-\pi_{i_j}$ up steps, and finally ends with $n+1-i_k$
down steps. It can be checked that this is indeed a bijection between $132$-avoiding permutations and Dyck paths.
The bijection we are looking for is $\bijb:=\krat^{-1}\circ\pre$.

Each young leaf of $T$ produces an occurrence of $UDU$ in $\pre(T)$. Such an occurrence appears in $\krat(\pi)$
for each pair of consecutive left-to-right minima. These two elements, together with the entry of $(n+1)\pi$ immediately to their left, form
a decreasing sequence of three consecutive elements (a double descent). To see that these are the only
double descents of $(n+1)\pi$, notice that from the structure of $132$-avoiding permutations
it follows that if $\pi_j>\pi_{j+1}$ is a descent of $\pi$, then $\pi_{j+1}$ must be a left-to-right minimum.

The reasoning for young leaves is similar. They correspond in
$\pre(T)$ to occurrences of $UDD$ and possibly a $UD$ at the end.
Equivalently, to occurrences of $UDD$ in $\pre(T)D$ (i.e., the Dyck
path $\pre(T)$ with a $D$ step appended at the end). Each of these
occurrences marks the start of a maximal sequence of at least two
consecutive $D$ steps in $\pre(T)D$, and each such sequence
corresponds to an ascending run of $\pi(n+1)$.
\end{proof}

For example, if $T$ is again the tree in Figure~\ref{fig:tree}, then
the corresponding $132$-avoiding permutation is
$\pi=\bijb(T)=(11,10,12,13,9,5,6,7,8,3,2,1,4)$. Note that
$(n+1)\pi=(14,\pi)$ has four double descents, namely $(14,11,10)$,
$(13,9,5)$, $(8,3,2)$ and $(3,2,1)$, and $(\pi,14)$ has three
ascending runs, namely $(10,12,13)$, $(5,6,7,8)$ and $(1,4,14)$.

\bs

There is another well-known bijection between plane trees and Dyck paths, which we denote $\dgr$. Given a tree $T$, traverse it in preorder (from
left to right) and build $\dgr(T)$ as follows.
For each node with $r$ children, draw $r$ up steps followed by one down step; except for the last leaf, for which we do
not draw anything. For example, the path corresponding to the tree in Figure~\ref{fig:tree} is $\dgr(T)=UUUUDUUUDDDDUDUDUDDDUDUUDD$.

Define a \emph{drop} of a Dyck path to be a maximal succession of at least two consecutive $D$ steps, and a \emph{triple fall}
to be an occurrence of $DDD$. Then the bijection
$\dgr$ maps each old leaf of $T$ to a drop of $\dgr(T)D$, and each young leaf to a triple fall of $\dgr(T)D$. In the example from the paragraph above,
$\dgr(T)D$ has three drops and four triple falls.

Following very similar arguments to the ones in Propositions \ref{prop:st321} and \ref{prop:st132}, but using the bijection $\dgr$ instead of $\pre$,
we obtain the next two results.

\begin{prop}\label{prop:dgr321}
There is a bijection $\bijc:\T_n\longrightarrow\S_n(321)$ such that, if $T\in\T_n$
and $\pi:=\bijc(T)\in\S_n(321)$, then
\ben
\item $\#$ of young leaves of $T = \#$ of pairs of consecutive deficiencies of $\pi$ ($+1$ if $\pi_n<n$),
\item $\#$ of old leaves of $T = \#$ of weak excedances of $\pi$ not followed by another weak excedance.
\een
\end{prop}

\begin{prop}
There is a bijection $\bijd:\T_n\longrightarrow\S_n(132)$ such that, if $T\in\T_n$
and $\pi:=\bijd(T)\in\S_n(132)$, then
\ben
\item $\#$ of young leaves of $T = \#$ of double ascents of $\pi(n+1)$,
\item $\#$ of old leaves of $T = \#$ of ascending runs of $\pi(n+1)$.
\een
\end{prop}

\section{Refinements of two combinatorial identities} \label{sec:coker}

In \cite{Cok} Coker established the following two identities, involving teh Narayana and the Catalan numbers:

\beq \label{cok1} \sum_{k=1}^n \frac{1}{n}\binom{n}{k}\binom{n}{k-1} 4^{n-k}=
\sum_{k=0}^{\lfloor(n-1)/2\rfloor} C_k \binom{n-1}{2k} 4^k 5^{n-2k-1}, \eeq
\beq \label{cok2} \sum_{k=1}^n \frac{1}{n}\binom{n}{k}\binom{n}{k-1} x^{2k}(1+x)^{2n-2k}=
x^2 \sum_{k=0}^{n-1} C_{k+1} \binom{n-1}{k} x^k (1+x)^k,
\eeq

He stated the open problem of finding a combinatorial interpretation of these identities. In
\cite{CYY}, Chen, Yan and Yang proved these identities combinatorially by applying 
our bijection $\bj$ to weighted plane trees.
In this section we use the properties of $\bj$ given in Proposition~\ref{prop:corresp} to obtain refinements
of the identities (\ref{cok1}) and (\ref{cok2}).

\begin{theorem} \label{th:cok1ref} For $n\ge1$, we have
\beq \label{cok1ref} \sum_{i=1}^n\sum_{j=0}^{n-2i+1} \frac{1}{n}\binom{n}{i}\binom{n-1}{j}\binom{n-i-j}{i-1} x^{i-1}y^{j}=
\sum_{k=0}^{\lfloor(n-1)/2\rfloor} C_k \binom{n-1}{2k} x^k (1+y)^{n-2k-1}. \eeq
\end{theorem}

\begin{proof}
We use a very similar reasoning to the one given in \cite{CYY} to prove equation (\ref{cok1}).
It will be convenient to use the term \emph{critical leaf} to denote
the last old leaf that we encounter when we traverse a plane tree in preorder.
Given a plane tree $T$ with $n$ edges, assign weights to the vertices of $T$ as follows:
young leaves are given weight $y$, old leaves other than the critical one are given weight $x$,
and the rest of the vertices (including the critical leaf) are given weight $1$. The weight of $T$ is the
product of the weights of its vertices. Then, the left hand side of (\ref{cok1ref}) is the sum of the weights
of all plane trees with $n$ edges.

By Proposition~\ref{prop:corresp}, $\bj$ is a weight preserving bijection between 
the set of weighted plane trees on $n$ edges, with weights given as above, and the set of weighted 2-Motzkin paths 
of length $n-1$ where weights are assigned as follows: $U$ steps are given weight $x$, 
$R$ steps are given weight $y$, and all the remaining steps are given weight $1$, defining the weight 
of a 2-Motzkin path to be the product of weights of its steps.
We claim that the right hand side of (\ref{cok1ref}) is the sum of the weights of all 2-Motzkin paths of length $n-1$.
Indeed, let $k\le\lfloor(n-1)/2\rfloor$ and consider the weighted 2-Motzkin paths with $k$ up steps and $k$ down steps.
These up and down steps from a 
Dyck path of length $2k$, and the positions of these $2k$ steps can be chosen in $\binom{n-1}{2k}$ ways. 
They contribute $x^k$ to the weight of the path.
The remaining $n-2k-1$ steps are either $R$ or $B$ steps. Since $R$ steps have weight $y$ and $B$ steps have 
weight $1$, the total contribution of the the horizontal steps in paths with $k$ up steps is $(1+y)^{n-2k-1}$. 
This justifies the right hand side.
\end{proof}

With the subsitution $y=x$ in equation (\ref{cok1ref}) we recover the result proved in \cite{CYY}, and 
the particular case $y=x=4$, together with the symmetry of the Narayana numbers, yields equation (\ref{cok1}).
A refinement of the second identity (\ref{cok2}) is given next.

\begin{theorem}
For $n\ge1$, we have
\beq \label{cok2ref} \sum_{i=1}^n\sum_{j=0}^{n-2i+1} \frac{1}{n}\binom{n}{i}\binom{n-1}{j}\binom{n-i-j}{i-1} 
x^{2(i-1)} y^j z^{n-2i-j+1} =
\sum_{k=0}^{n-1} C_{k+1} \binom{n-1}{k} x^k (y+z-2x)^{n-1-k}. \eeq
\end{theorem}

\begin{proof}
Again we apply the same ideas used in \cite{CYY} to prove equation (\ref{cok2}).
Recall the definition of critical leaf from the proof or Theorem~\ref{th:cok1ref}.
Given a plane tree $T$ with $n$ edges, assign weights to the vertices of $T$ in the following way.
Old leaves other than the critical one are given weight $x$, the parents of such leaves are given weight $x$
as well, young leaves are given weight $y$, the critical leaf and its parent are given weight $1$, and the rest of the 
vertices are given weight $z$. As before, the weight of $T$ is the
product of the weights of its vertices. Notice that two different old leaves cannot have the same parent, so
the weight of a tree with $i$ old leaves and $j$ young leaves is $x^{2(i-1)} y^j z^{n-2i-j+1}$.
The left hand side of (\ref{cok2ref}) is the sum of the weights of all plane trees with $n$ edges.

By Proposition~\ref{prop:corresp}, a tree with $i$ old leaves and $j$ joung leaves is mapped by $\bj$ to a 2-Motzkin
path with $i-1$ up steps, $i-1$ down steps, $j$ horizontal $R$ steps, and $n-2i-j+1$ horizontal $B$ steps.
To make $\bj$ a weight preserving bijection between plane trees on $n$ edges with the above weights 
and 2-Motzkin paths of length $n-1$, we assign weights to the steps of a 2-Motzkin path by giving  weight $x$
to $U$ and $D$ steps, weight $y$ to $R$ steps, and weight $z$ to $B$ steps.

Consider now the set of 3-Motzkin paths of length $n-1$, where horizontal steps can be either red, blue or green
(call them $R$, $B$ and $G$ steps respectively). Assign weights to the steps by giving weight $y+z-2x$ to $G$ steps 
and weight $x$ to all the other steps.
This weight assignment in 3-Motzkin pahts has the property that the sum of the weights of an $R$ step, 
a $B$ step and a $G$ step equals the sum of the weights of an $R$ step and a $B$ step
in the assignment for 2-Motkin paths above (namely $x+y$), and also that $U$ and $D$ steps have the same 
weight $x$ in both assignments.
This implies that the sum of weights over all 2-Motkin paths with the above weight assignment equals the sum
of weights over all 3-Motzkin paths with this new assignment.
Therefore, all that remains is to show that the right hand side of (\ref{cok2ref}) is the total sum of the weights 
of 3-Motzkin paths of length $n-1$. But this is clear because if we fix the number of $G$ steps of a 3-Motzkin path
to be $n-1-k$, then the positions of these $G$ steps can be chosen in $\binom{n-1}{k}$ ways. The remaining steps,
$U$, $D$, $R$ and $B$, form a 2-Motkin path of length $k$, and the number of such paths is $C_{k+1}$.
\end{proof}

To recover identity (\ref{cok2}) we only need to substitute $x(1+x)$ by $x$, $x^2$ by $y$, and $(1+x)^2$ by $z$ in
equation (\ref{cok2ref}).

\section*{Acknowledgements}

We are grateful to Laura Yang for her valuable suggestions.
The work of W. Y. C. Chen was partially supported by the National Science 
Foundation, the Ministry of Education, and the
Ministry of Science and Technology of China.
The work of S. Elizalde was partially supported by the Ministry of Foreign Affairs of Spain and the AECI.

\end{document}